\newtheorem{theorem}{\sc Theorem}

\newtheorem{claim}{\sc Claim}[theorem]
 
\newcommand{\proof}{\noindent {\sc Proof. }}

\def\rest{\mathord{\restriction}}
\newcommand{\force}{\Vdash}
\newcommand{\open}{\Bbb}

\newcommand{\dom}{\mbox{\rm dom}}

\newcommand{\se}{\subseteq}
\newcommand{\set}[2]{\{#1 \colon #2\}} 
\newcommand{\fin}{$\Box$\par\medskip} 

\renewcommand{\to}{\rightarrow}
\def\rge{\mathop{\rm rge}}

\newcommand{\ga}{\alpha}
\newcommand{\gb}{\beta}

\newcommand{\gd}{\delta}

\newcommand{\gk}{\kappa}

\newcommand{\gt}{\tau}

\newcommand{\go}{\omega}
 
\newcommand{\gF}{\Phi}

\newcommand{\ha}{\aleph}

\newcommand{\oP}{{\open P}}

\newcommand{\tb}{{\tilde b}}
\newcommand{\tr}{{\tilde r}}
\newcommand{\tS}{{\tilde S}}
\newcommand{\fA}{{\frak A}}
\newcommand{\fB}{{\frak B}}
\newcommand{\frC}{{\frak C}}

\newcommand{\gogo}{{}^{<\go_1}\go_1}
\newcommand{\seq}[2]{\langle #1\colon #2\rangle}
\newcommand{\EF}{Ehrenfeucht-Fra\"{\i}ss\'e}

\newcommand{\hk}{({\rm H}(\gk), {\in}, <^*)}

\documentstyle[12pt,amsfonts,amssymb]{article}

\title{The Canary Tree}

\author{Alan H. Mekler\thanks{Research partially supported by NSERC
grant A8948. Research on
this paper was begun while both authors were visiting
MSRI.}\\Department of Mathematics and Statistics\\ Simon Fraser 
University\\Burnaby, B.C., V5A 1S6\\ Canada \and Saharon
Shelah\thanks{Publication
\#398. Research supported by the BSF and NSF.}\\Institute of
Mathematics\\The Hebrew University\\Jerusalem 91904, Israel\\and\\Department
of Mathematics\\ Rutgers University\\New Brunswick, New Jersey
08903\\USA}
\date{}

\begin{document}
\maketitle
\begin{abstract}
A {\em canary tree} is a tree of cardinality the continuum which has no
uncountable branch, but  gains a branch whenever a stationary set is
destroyed (without adding reals). Canary trees are important in
infinitary model theory. The existence of a canary tree is
independent of ZFC $+$ GCH.

\medskip\noindent
{\em 1991 Mathematics Subject Classification.} Primary 03E35;
Secondary 03C75.
\end{abstract}

	A canary tree is a tree of cardinality $2^{\ha_0}$ which
detects the destruction of stationary sets. (A stationary set is {\em
destroyed} in an extension if it is non-stationary in the extension.)
More exactly, $T$ is a {\em canary tree} if $|T| = 2^{\ha_0}$, $T$ has
no uncountable branch, and in any extension of the universe in which no new
reals are added and in which some stationary subset of $\go_1$ is
destroyed, $T$ has an uncountable
branch. (We will give an equivalent characterization below
which does not mention extensions of the universe.) The existence of
a canary tree is most interesting under the assumption of CH (if $2^{\ha_0}
= 2^{\ha_1}$ it is easy to see, as we will point out, that there is a
canary  tree.) The existence or non-existence of a canary tree has
implications for the model theory of structures of cardinality
$\ha_1$ and for the descriptive set theory of ${}^{\go_1}\go_1$
(\cite{MV}). The canary tree is named after the miner's canary.

	In this paper, we will explain the significance of the
existence of a canary tree in model theory and prove that the
existence of a canary tree is independent from ZFC $+$ CH.

	As is well known the standard way to destroy a stationary
costationary subset of $\go_1$ is to force a club through its
complement using as conditions closed subsets of the complement
(\cite{BHK}). More precisely if $S$ is a stationary subset of $\go_1$
we can define $T_S = \set{C}{C \mbox{ a closed countable subset of }
S}$, where the order is end-extension. If $S$ is costationary then
$T_S$ has no uncountable  branch but when we force with $T_S$ we add
no reals but do add a branch through $T_S$. Such a branch is a club
subset of $\go_1$ which is contained in $S$. (In \cite{BHK}, the
forcing to destroy a stationary costationary set $E$ is exactly
$T_{\go_1 \setminus E}$.) Notice that $T_S$ detects the destruction of
$\go_1 \setminus S$ in the sense that in any extension of the universe
with no new reals and in which  $\go_1 \setminus S$ is non-stationary,
$T_S$ has a branch. 

	These elementary observations imply that if $2^{\ha_0}
= 2^{\ha_1}$ then there is a canary tree. The tree can be constructed
by having disjoint copies of $T_S$ sitting above a common root where
$S$ ranges through the stationary costationary subsets of $\go_1$. In
fact any canary tree must almost contain the union of all the $T_S$,
in the following weak sense. 

\begin{theorem}
\label{equiv}
Suppose $T$ is a tree of $2^{\ha_0}$ with no uncountable branch. Then
$T$ is a  canary tree if and only if for any  stationary
costationary set $S$
there exists a sequence $\seq{X_\ga}{\ga < \go_1}$ of maximal
antichains of $T_S$  and there is an order-preserving function $f\colon
\bigcup_{\ga< \go_1} X_\ga \to T$. Furthermore  $\seq{X_\ga}{\ga <
\go_1}$ and $f$ are such that: if $\ga < \gb$ and $s
\in X_\ga$, $t \in X_\gb$ then either $s$ and $t$ are incomparable or
$s < t$; if $\gd$ is a limit ordinal and $t \in X_\gd$ then $t =
\sup \{s < t \colon s \in X_\gb, \gb < \gd\}$; and $f$ is continuous.
(Note that these conditions imply that for all $u \in T_S$ there is
$\gd$ and $t\in X_\gd$ such that $u < t$.)

\end{theorem}

\proof First assume that for every stationary costationary set there
is such a sequence of  antichains and such a function. Suppose that
$E$ is a stationary set which is destroyed in an extension of the
universe with no new reals. Let $S = \go_1 \setminus E$ and let $f$
and $\seq{X_\ga}{\ga < \go_1}$ be as guaranteed.  Let $C$ be a club in
the extension which is contained in $S$. Choose an increasing sequence
$\seq{s_\ga}{\ga < \go_1}$ of elements of $T_S$ so that for all $\ga$,
$s_{\ga +1}$ is greater than some member of $X_\ga$ and $\max s_\ga
\in C$. The choice of such a sequence is by induction. There is no
problem at successor steps. At a limit ordinal $\gd$, we can continue
since $\sup \bigcup_{\ga < \gd} s_\ga \in C$ and hence in $S$. Also
since no new reals are added $\bigcup_{\ga < \gd} s_\ga \cup\sup
\bigcup_{\ga < \gd} s_\ga \in T_S$. Let $b$ be the uncountable branch
through $T_S$ determined by $\seq{s_\ga}{\ga < \go_1}$. So in the
extension $f\mbox{''}(b \cap \bigcup_{\ga< \go_1} X_\ga)$ is an
increasing uncountable subset of $T$.

	Now suppose that $T$ is a canary tree. 
	Let $S$ be a stationary costationary set. Since forcing with
$T_S$ destroys a stationary set there is $\tb$ a $T_S$-name for a
branch of $T$. We will inductively define the  sequence $\seq{X_\ga}{\ga
< \go_1}$ of maximal antichains of $T_S$. Let $0$ denote be the root of
$T$. Define $X_0 = \{0\}$. In general, let $Y_\ga = T \setminus
\bigcup_{\gb < \ga} X_\gb$ and let  $D_\ga = \set{t \in Y_\ga}{t \mbox{
decides } \tb\rest\ga}$. Let $X_\ga$ be the set of minimal elements of
$D_\ga$. Since $D_\ga$ is dense, $X_\ga$ is a maximal antichain.
For $t \in X_\ga$, choose $s$ so that $t \force s = \tb\rest\ga$ and
let $f(t) = s$. \fin  

	It is possible to improve the theorem above to show that
$T$ is a canary tree if and only if for every stationary costationary set $S$
there is an order preserving function from $T_S$ to $T$ (\cite{MV}). In
fact  when we show that it is consistent with GCH that there is a
canary tree $T$, we will construct  for every stationary costationary set $S$
an order preserving function from $T_S$ to $T$. It is also
worth noting that we get an equivalent definition if we only demand
that a canary tree have cardinality at most 
$2^{\ha_0}$, since if $T$ is a tree of cardinality less
than $2^{\ha_0}$, then forcing with $T_S$ adds no new branch to $T$.

\section{The Canary Tree and \EF\ Games}

	A central idea in the Helsinki school's approach to finding an
analogy at $\go_1$ of  the theory of $L_{\infty \go}$ is the notion
of an \EF\ game of length $\go_1$ (see \cite{HV} for more details and
further references). Given two models, $\fA$ and $\fB$,
two players, an isomorphism player and a non-isomorphism player,
alternately choose elements from $\fA$ and $\fB$. In its primal form
the game lasts $\go_1$ moves and the isomorphism player wins if an
isomorphism between the chosen substructures has been constructed. The
analogue of Scott's theorem is the trivial result that two structures
of cardinality $\ha_1$ are isomorphic if and only if the isomorphism
player has a winning strategy. In the search for an analogue of Scott
height, trees with no uncountable branches play the role of ordinals.
More exactly suppose that $T$ is a tree and $\fA$ and $\fB$ are
structures. The game ${\cal G}_T(\fA, \fB)$ is defined as follows. At
any stage the non-isomorphism player chooses an element from either
$\fA$ or $\fB$ and a node of $T$ which lies above the nodes this
player has already chosen. The isomorphism player replies with an
element of $\fB$ if the non-isomorphism player has played an element
of $\fA$ and an element of $\fA$ if the non-isomorphism player has
played an element of $\fB$. In either case the move must be such that
the resulting sequence of moves from $\fA$ and $\fB$ form a partial
isomorphism. The first player who is unable to move loses. In analogy
with Scott height if $\fA$ and $\fB$ are non-isomorphic structures of
cardinality $\ha_1$ then there is a tree of cardinality at most $2^{\ha_0}$
with no uncountable branches such that the non-isomorphism player has
a winning strategy in ${\cal G}_T(\fA, \fB)$. (The tree $T$ can be
chosen to be minimal.) A defect in the analogy with Scott height is
that the choice of the tree depends on the pair $\fA, \fB$ and cannot
in general be chosen for  $\fA$ to work for all $\fB$ (\cite{HT}). 

\noindent {\sc Definition.} Suppose $\fA$ is a structure of
cardinality $\ha_1$. A tree $T$ is called a {\em universal
non-equivalence tree} for $\fA$ if $T$ has no uncountable branch and
for every non-isomorphic $\fB$ of cardinality $\ha_1$ the
non-isomorphism player has 
a winning strategy in   ${\cal G}_T(\fA, \fB)$.

\smallskip

	As we have mentioned there are structures 
for which there is no universal non-equivalence tree of cardinality
$\ha_1$. However for some natural structures such as free groups (or
free abelian groups)
or $\go_1$-like dense linear orders the existence of a universal
non-equivalence tree of cardinality $2^{\ha_0}$ is equivalent to the
existence of a canary tree. We will only explain the case of
$\go_1$-like dense linear orders, the case of groups is
similar. 

	Recall the classification of $\go_1$-like dense linear orders
with a left endpoint. Let $\eta$ represent the rational order type and
for $S \se \go_1$ let $\gF(S) = 1 + \eta + \sum_{\ga < \go_1}
\gt_\ga$, where $\gt_\ga = 1 + \eta$ if $\ga \in S$ and $\gt_\ga =
\eta$ otherwise. It is known that any $\go_1$-like dense linear order
is isomorphic to some $\gF(S)$ and that for $E, S \se \go_1$, $\gF(S)
\cong \gF(E)$ if and only if the symmetric difference of $E$ and $S$
is nonstationary.

\begin{theorem}
There is a universal non-equivalence tree of cardinality $2^{\ha_0}$ for
$\gF(\emptyset)$ if and only if there is a canary tree. 
\end{theorem}

\proof Assume  that $T$ is a
universal non-equivalence tree of cardinality $2^{\ha_0}$ for
$\gF(\emptyset)$. Consider $E$, a stationary  costationary set.
Work now in an extension of the universe in which $E$ is
non-stationary and there are no new reals. In that universe, $\gF(E)
\cong \gF(\emptyset)$. In that universe the isomorphism player can
play the isomorphism against the winning strategy of the
non-isomorphism player in ${\cal G}_T(\gF(\emptyset),\gF(E))$. At
each stage, both players will have a move. So the game will last
$\go_1$ moves and the non-isomorphism player will have chosen an
uncountable branch through $T$. Hence $T$ is a canary tree.

	Now suppose that $T$ is a canary tree. Let $T' = T+ 2$ (i.e.,
a chain of length 2 is added to the end of every maximal branch of
$T$). We claim that $T'$ is a universal non-equivalence tree for
$\gF(\emptyset)$. Suppose $E$ is a stationary set. The case where $E$
is in the club filter is an easier version of the following argument.
Assume that $S$ is stationary where $S = \go_1 \setminus E$.
To fix notation let $\gF(\emptyset) = 1 + \eta + \sum \gt_\ga$ and
$\gF(E) = 1 + \eta + \sum \mu_\ga$.  Let $\seq{X_\ga}{\ga < \go_1}$
and $f\colon T_S \to T$ be as in Theorem~\ref{equiv}. Let $X =
\bigcup_{\ga < \go_1} X_\ga$. The winning strategy for the
non-isomorphism player consists of choosing an increasing sequence
$s_\ga \in X$, playing $f(s_\ga)$ as the move in the tree $T$ and
guaranteeing that at every limit ordinal $\gd$ if $A$ is the subset
of $\gF(\emptyset)$ which has been played (by either player) and $B$
is the subset of $\gF(E)$ which has been played then $\sup\bigcup_{\ga <
\gd} s_\ga = \sup \set{\gb}{a \in \gt_\gb, a \in A} = \sup \set{\gb}{b
\in \mu_\gb, b \in B}$. The non-isomorphism player continues this way
as long as possible. When there are no more moves following this
recipe $\sup\bigcup_{\ga < 
\gd} s_\ga$ is an ordinal in $E$. In that case $B$ has a least upper
bound but $A$ doesn't. So the non-isomorphism player only needs two
more moves to win the game. \fin

	The above argument also shows that if there is a canary tree
then the $\go_1$-like dense linear orders share a 
universal non-equivalence tree of cardinality $2^{\ha_0}$. 

\section{Independence Results}

\begin{theorem}
\label{notree}
    It is consistent with GCH that there is no canary tree.
\end{theorem}

\proof Begin with a model of GCH and add $\ha_2$ Cohen subsets to
$\go_1$. In the extension GCH continues to hold. Suppose $T$ is a tree
of cardinality $\ha_1$ which has no uncountable branch. Since the
forcing to add $\ha_2$ Cohen subsets of $\go_1$ satisfies the
$\ha_2$-c.c., $T$ belongs to the extension of the universe by $\ha_1$
of the subsets. By first adding all but one of the subsets we can work
in $V[X]$ where $X$ is a Cohen subset of $\go_1$ and $T$ is in $V$.
Note that $X$ is a stationary costationary subset of $\go_1$. Let
$\oP$ be the forcing for adding a Cohen generic subset of $\go_1$ and
let $Q$ be the $\oP$-name for $T_X$. It is easy to see that $\oP \ast
Q$ is essentially $\go_1$-closed. Hence forcing with $\oP \ast Q$
doesn't add a branch through $T$. So neither does forcing with $T_X$
over $V[X]$. But forcing with $T_X$ destroys a stationary set, namely,
$\go_1 \setminus X$. \fin

	It remains to prove the consistency of GCH together with the
existence of a canary tree. The proof has two main steps, we first
force a very large subtree of $\gogo$. At limit ordinals we will
forbid at most one branch from extending. Having created the tree we
will then iteratively force order preserving maps of  $T_S$ into the
tree as $S$ varies over all stationary costationary sets.

\begin{theorem}
    It is consistent with GCH that there is a canary tree.
\end{theorem}

\proof Assume that GCH holds in the ground model. To begin define 
$Q_0$ to be 
$$\set{f\colon \lim(\go_1) \to \gogo}{\dom f \mbox{ is countable and
for all }\gd \in \dom(f), 
f(\gd) \in {}^\gd\gd}.$$
If $G_0$ is $Q_0$-generic, we can identify $G_0$ with $\bigcup_{f \in
G_0} \rge f$. Let ${\frC} = \set{s \in \gogo}{\mbox{for
all }\gd \leq \ell(s), s\rest\gd \notin G_0}$. It is easy to see that
in $V[G_0]$, $\frC$ has no uncountable branch and that $V[G_0]$
has no new reals. (In fact forcing with $Q_0$ is the same as adding a
Cohen subset of $\go_1$, so the claims above follow.)

	To complete the proof we need to force embeddings of $T_S$
into $\frC$ as $S$ ranges over stationary sets. Suppose that we are in
an extension of the universe which includes a generic set for $Q_0$
and has no new reals. Fix a stationary set
 $S$. An element $t$ of $\frC$ is called an $S$-{\em node} if for
every limit ordinal $\ga \notin S$, if $\ga \leq \ell(t)$ then  $t\rest
\ga \notin {}^\ga\ga$.  Notice that any $S$-node has successors of
arbitrary height, since if $s$ is an $S$-node of height $\ga$ and
$\gd$ is a limit ordinal greater than $\ga$, then any extension of
$s^\frown\langle\gd\rangle$ of length at most $\gd$ is an $S$-node.
The poset $\oP(S)$ will consist of pairs $(g, X)$ where $X$ is a
countable subset of $\gogo$ such that each element of $X$ is of
successor length and $g$ is a partial order preserving map from $T_S$
to the $S$-nodes of $\frC$ whose domain is a countable subtree of
$T_S$. Further $(g, X)$ has the following properties. 
\begin{enumerate}
    \item if $c \in \dom(g)$ and $t \in X$ then $t \not\se g(c)$
    \item if $c_0 < c_1 < \ldots$ is an increasing sequence of
elements of $\dom(g)$ then $\bigcup_{n<\go} g(c_n) \in \frC$.
\end{enumerate}

	If $(g, X)$ is a condition let $o(g, X)$ be the
$\sup\set{\ell(t)}{t \in X \mbox{ or } t \in \rge(g)}$. A condition $(h, Y)$
extends $(g, X)$ if 
\begin{enumerate}
    \item $g \se h$,
    \item if $c \in \dom(h) \setminus \dom(g)$, then $\ell(h(c)) >
o(g, X)$,
    \item $X \se Y$.
\end{enumerate}

\begin{claim}
\label{proper}
    The poset $\oP(S)$ is proper. 
\end{claim}

	Suppose $\gk$ is some suitably large cardinal, $N \prec \hk$,
where $<^*$ is a well-ordering of the model, $N$ is countable, and
$\oP(S) \in N$. We need to show that for every $p \in N \cap \oP(S)$
there is an $N$-generic extension. Let $\gd = N \cap \go_1$. Let $f$
be the $Q_0$-generic function and $t = f(\gd)$. There are two cases to
consider. Either there is a successor ordinal $\ga < \gd$ so that $\ga
> o(p)$ and  $t \rest \ga \in N$ or not. Let $p = (g, X)$. If such an
ordinal $\ga$ exists let $p_{-1} = (g, X \cup \{t\rest\ga\})$,
otherwise let $p_{-1} = p$. Now define a sequence $p_{-1}, p_0,
\ldots, p_n \ldots$ of increasingly stronger conditions so that (for
$n \geq 0$) $p_n$
is in the $n^{\rm th}$ dense subset of $\oP(S)$ which is an element of
$N$. Let $p_n = (g_n, X_n)$ and $q = (h, Y)$ where $h = \bigcup_{n <
\go} g_n$ and $Y =  \bigcup_{n < \go} X_n$. To finish the proof it
suffices to see that $q \in \oP(S)$. The only point that needs to be
checked is to verify that if $c_0 < c_1 < \ldots \in \dom(h)$ then
$\bigcup_{n < \go} h(c_n) \in \frC$. If there is $m$ such that $c_n
\in \dom(g_m)$ for all $n$, then we are done. Otherwise, by the second
property of being an extension, $\sup\set{\ell(h(c_n)}{n < \go} \geq
\sup\set{o(g_m, X_m)}{m < \go}$. However for all $\ga < \gd$ there is
a dense set $D$ such that $(g, X) \in D$ implies $o(g, X) > \ga$. As
$D$ is definable using parameters from $N$, $D \in N$. Furthermore
since the sequence of conditions meets every dense set in $N$,
$\sup\set{o(g_m, X_m)}{m < \go} \geq\gd$. Finally each $h(c_n) \in N$,
so $\ell(h(c_n)) < \gd$ for all $n$. These facts give the equation,
$\sup\set{\ell(h(c_n)}{n < \go} = \gd$. (In the remainder of
the paper we will try to point out where a density argument is needed
but we will not give it in such detail.) 
By the choice of $p_{-1}$ and the property 1 of the definition
of $\oP(S)$, $t \neq \bigcup_{n < \go} h(c_n)$. \fin 

	Our forcing will be an iteration with countable support of
length $\go_2$. As usual we will let $\oP_i$ be the forcing up to
stage $i$ and will force with $Q_i$, a $\oP_i$-name for a poset. We
have already defined $Q_0$. For $i$ greater than $0$, we take $\tS_i$
a $\oP_i$-name for a stationary costationary set and let $Q_i$ be the
$\oP_i$-name for $\oP(\tS_i)$. By  Claim~\ref{nnr},
forcing with $\oP_{\go_2}$ adds no reals. Also since each $Q_i$ is
forced to have cardinality $\go_1$, if we enumerate the $\tS_i$
properly every  stationary costationary set in the final forcing
extension will occur as the interpretation of some $\tS_i$.

\begin{claim}
\label{nnr}
For all $i \leq \go_2$, forcing with $\oP_i$ adds no new reals.
\end{claim}

	The proof is by induction on $i$. The case $i = 1$ is easy.
For successor ordinals the proof can be done along the same lines as
Claim~\ref{proper}, or by a modification of the limit ordinal case which
we do below. Suppose now that $i$ is a limit ordinal and $\tr$ is a
$\oP_i$-name for a real. Consider any condition $p$. We must show that
$p$ has an extension which determines all the values of $\tr$. Choose
a countable $N$ so that $N \prec \hk$ and $p, \oP_i, \tr \in N$.
Let $ p = p_{-1}, p_0, p_1, \ldots$ be a sequence of increasingly
stronger conditions in $N$ so that $p_n$ is in the $n^{\rm th}$ dense
subset of $\oP_i$ which is an element of $N$.  Let $\gd = N \cap
\go_1$.   There is an obvious upper
bound $q$ for the sequence. Of course $q$ is not a condition.  We
would like to extend $q$ to a condition $q'$ by choosing some $t \in
{}^\gd\gd$, letting $q'(0) = 
q(0)^\frown\langle\gd, t\rangle$ and letting $q'(i) = q(i)$ for $i >
0$. Choose $t\in {}^\gd\gd$ so that $t\rest\go \notin N$. By a
density argument we can show that for all $i$ and $n$, if $p_n\rest
i\force c \in \dom p$, then there are $m$, $g$, $X$ and $s \in
{}^\go\go \cup N$ so that $p_m \force p_n(i) = (g, X) \mbox{ and }
g(c) = s$. It is straightforward to see that $t$ is as desired. (See
the proof of Claim~\ref{nobranch} for a similar but more detailed
argument.) \fin

	Let $G_{\go_2}$ be $\oP_{\go_2}$-generic. We have shown that in
$V[G_{\go_2}]$, for every stationary set $S$ there is an order
preserving map from $T_S$ to $\frC$. To finish the proof we must
establish the following claim.

\begin{claim}
\label{nobranch}
   In $V[G_{\go_2}]$, $\frC$ has no uncountable branch.
\end{claim}

	 Suppose that $\tb$ is forced (for simplicity) by the empty
condition to be an uncountable branch of $\gogo$. We will show that
there is a dense set of conditions which forces that $\tb$ is not a
branch of $\frC$. Hence $\frC$ has no uncountable branch. Fix a
condition  $p \in \oP$.  Choose a countable $N$ so that
$N \prec \hk$ and $p, \oP_{\go_2}, \tb \in N$. Let $ p = p_{-1}, p_0, p_1,
\ldots$ be a sequence of increasingly stronger conditions in $N$ so that
$p_n$ is in the $n^{\rm th}$ dense subset of $\oP_{\go_2}$ which is an
element of $N$.  Let $\gd = N \cap \go_1$. The sequence $\seq{p_n}{n <
\go}$ determines a value for $\tb\rest \gd$. Let this value be $t$.
There is an obvious upper bound $q$ for the sequence. Of course $q$ is
not a condition.  We would like to extend $q$ to a condition $q'$ by
letting $q'(0) = q(0)^\frown\langle\gd, t\rangle$ and letting $q'(i) =
q(i)$ for $i > 0$. We will show by induction on $i$ that $q'\rest i$
is a condition in $\oP_i$.

	The case $i = 1$ and limit cases are easy. So we can assume
that $i \in N$ and $q'\rest i \in \oP_i$. Since forcing with $\oP_i$
adds no new reals and $N$ is an elementary submodel of $\hk$, for all
$n$ there is $m$ and $(g_n, X_n)$ so that $p_m\rest i \force p_n(i) =
(g_n, X_n)$. Hence $q'\rest i \force q'(i) = (h, Y)$, where $h
=\bigcup_{n< \go} g_n$ and $Y =  \bigcup_{n < \go} X_n$.
Suppose now that $c_0 < c_1 < \ldots \in \dom(h)$. We need to show
that $q'\rest i\force \bigcup_{n < \go} h(c_n) \in \frC$. If there is some
$m$ so that $c_n \in \dom(g_m)$ for all $n$, then we are done as in
Claim~\ref{proper}. Otherwise $\ell(\bigcup_{n < \go} h(c_n)) = \gd$
and we only need to show that $\bigcup_{n < \go} h(c_n) \neq t$. 

	Notice that for all $\ga < \gd$, $q'\rest i\force (\bigcup_{n < \go}
h(c_n))\rest \ga \mbox{ is an }\tS_i\mbox{-node}$. We
will show that there is $\ga < \gd$ so that then $q'\rest i
\force t\rest\ga \mbox{ is not an $\tS_i$-node}$. This will complete the
proof.

	Let $G = \set{p \in N \cap \oP_{\go_2}}{ \mbox{there is $n$ so
that $p_n$ extends }p}$. By the choice of the sequence, $G$ is
$N$-generic. Note that by Claim~\ref{proper} and the iteration lemma
for proper forcing (or by a direct
argument similar to Claim~\ref{proper}), $\oP_{\go_2} \force
\tS_i \mbox{ is costationary}$. Hence for all $i \in N$, $N[G] \models
\tS_i^G \mbox{ is 
costationary}$ and $N[G] \models \set{\ga}{\tb^G\rest\ga \in
{}^\ga\ga} \mbox{ is a club}$. Hence 
$$N[G] \models \mbox{there is a
limit ordinal } \ga \mbox{ so that } \tb^G\rest\ga \in {}^\ga\ga \mbox{ and }
\ga\notin \tS_i^G.$$ 
By the forcing theorem there is some $n$ so that $p_n\rest i \force
t\rest\ga\in {}^\ga\ga \mbox{ and } \ga\notin \tS_i$. So we have
shown  $q'\rest i\force t\rest\ga \mbox{ is not an $\tS_i$-node}$,
which was our goal. 
\fin

	Note in the proof above it was necessary to force the
embeddings. The forcing $Q_0$ is the same as adding a Cohen subset of
$\go_1$. So if we add two Cohen subsets of $\go_1$ and use one to
construct the tree, then, by the proof of Theorem~\ref{notree} the
other one gives a stationary set which can be destroyed without adding an
uncountable branch.


\begin{thebibliography}{99}
\bibitem{BHK} Baumgartner, J., Harrington, L. and Kleinberg, G. {\em
Adding a closed unbounded set}, J. Symbolic Logic {\bf 41}(1976) 481--482.

\bibitem{HT}  Hyttinen, T. and Tuuri, H. {\em Constructing strongly
equivalent nonisomorphic models for unstable theories}, Ann. Pure and
Appl. Logic {\bf 52}(1991) 203--248.

\bibitem{HV} Hyttinen, T and V\"a\"an\"anen, J. {\em On Scott and
Karp trees of uncountable models}, J. Symbolic Logic {\bf 55}(1990)
897--908. 

\bibitem{MV} Mekler, A. and V\"a\"an\"anen, J., {\em Trees and
$\Pi^1_1$-subsets of ${ }^{\omega_1}\omega_1$}, submitted.
\end{thebibliography}
\end{document}